\documentclass{amsproc}
\usepackage{amssymb}
\usepackage[cmtip,all]{xy}
\usepackage[utf8]{inputenc}
\usepackage{longtable}
\usepackage{cite}
\usepackage[all,cmtip]{xy}
\usepackage{url}
\usepackage{hyperref}
\usepackage{cleveref}
\usepackage{setspace}
\usepackage{makecell}
\usepackage{float}
\usepackage{adjustbox}
\newcommand{\curvepage}[1]{\href{https://beta.lmfdb.org/ModularCurve/Q/#1/}{\texttt{#1}}}

\newtheorem{theorem}{Theorem}[section]
\newtheorem{Lemma}[theorem]{Lemma}
\theoremstyle{definition}

\theoremstyle{remark}
\newtheorem{remark}[theorem]{Remark}
\newcommand{\ConditionName}{LMFDB Conditions}
\newtheorem*{LMFDBCondition*}{LMFDB Conditions}
\newtheorem{question}[theorem]{Question}
\newtheorem{proposition}[theorem]{Proposition}
\newtheorem{conjecture}[theorem]{Conjecture}
\newtheorem{corollary}[theorem]{Corollary}
\numberwithin{equation}{section}

\begin{document}

\title[Quadratic Chabauty Experiments]{Quadratic Chabauty Experiments on Genus 2 Bielliptic Modular Curves in the LMFDB}
\author{Kate Finnerty}
\address{Kate Finnerty}
\curraddr{Boston University \\ 665 Commonwealth Avenue \\ Boston, MA 02215}
\email{ksfinn@bu.edu}

\begin{abstract}
We present results of quadratic Chabauty experiments on genus 2 bielliptic modular curves of Jacobian rank 2 that have recently been added to the LMFDB. We apply quadratic Chabauty methods over both the rationals and quadratic imaginary fields. In a number of cases, the experiments yielded algebraic irrational points among the set of mock rational points. We highlight specific notable examples, including the non-split Cartan modular curve $X_{ns}^+(15)$. Lastly, we offer a conjecture relating the level of the modular curve to the potential number fields over which points can arise
\end{abstract}

\maketitle

\section{Introduction}
\subsection{Background and Motivation}
For a projective, smooth, absolutely integral curve $X$ over $\mathbb{Q}$, let $X(\mathbb{Q})$ denote the set of its rational points. If $X$ is of genus 0, then $X(\mathbb{Q})$ is empty or infinite (and $X$ is isomorphic to $\mathbb{P}^1$). If the genus $g$ of $X$ is 1 and $X(\mathbb{Q})$ is nonempty, then $X$ is an elliptic curve, and Mordell's theorem states that $X(\mathbb{Q})$ is a finitely generated group \cite{mordell1922rational}. When the genus is 2 or greater, Faltings' theorem assures us that $X(\mathbb{Q})$ is finite \cite{faltings1983proof}. However, the proof does not provide an effective construction to \textit{find} these points. It is an ongoing question to provably compute $X(\mathbb{Q})$. 

A classical method to compute $X(\mathbb{Q})$ comes from the work of Chabauty \cite{chabauty1941points} and Coleman \cite{coleman1985effective}. Given a prime $p$ of good reduction for $X$ and a known $P\in X(\mathbb{Q})$, if the rank of the Jacobian variety $J$ of $X$ is less than the genus of the curve, we can  construct a finite set $X(\mathbb{Q}_p)_1$ that contains $X(\mathbb{Q})$. In practice, one can often succeed in provably computing all of $X(\mathbb{Q})$ by checking to see which $p$-adic points easily lift to rational points and then using a tool such as the Mordell-Weil sieve to rule out the remaining points from being rational.

In the case where the rank of the curve is not less than the genus, this fails to provide finite bounds. Kim developed a nonabelian Chabauty program using Selmer varieties and their iterated $p$-adic integrals \cite{kim2009Selmer}. Just as the classical method creates the set $X(\mathbb{Q}_p)_1$, nonabelian Chabauty produces analogous sets:
\begin{equation*}
    X(\mathbb{Q}_p)_1\supset X(\mathbb{Q}_p)_2\supset X(\mathbb{Q}_p)_3\supset\cdots\supset X(\mathbb{Q}_p)_n\supset X(\mathbb{Q}).
\end{equation*}

Kim conjectured that for $n$ sufficiently large, $X(\mathbb{Q}_p)_n$ is finite. He further conjectured that for sufficiently large $n$, we achieve the equality $X(\mathbb{Q}_p)_n=X(\mathbb{Q})$.

This work has led to a wave of development already in the first nonabelian case, $n=2$, known as quadratic Chabauty. The quadratic Chabauty set $X(\mathbb{Q}_p)_2$ is known to be finite when the rank of the Jacobian of the curve is less than $g-1+rk(NS(J))$, where $NS(J)$ denotes the N\'eron-Severi group of the Jacobian over $\mathbb{Q}$. Balakrishnan and Dogra establish this bound for finiteness in  \cite{Balakrishnan2018QC1} and together with M\"uller compute some explicit examples. Developments have led to computational tools for new and more complicated examples in recent years. In \cite{bianchi2020bielliptic}, Bianchi further developed quadratic Chabauty by using the $p$-adic sigma function instead of the double Coleman integral to compute the quadratic Chabauty locus for genus 2 bielliptic curves. Balakrishnan, Besser, Bianchi, and M{\"u}ller later generalized this to rational points on genus 2 bielliptic curves over number fields \cite{balakrishnan2021numberfields}. Bianchi and Padurariu then implemented and scaled these methods to compute the set of rational points on bielliptic curves of Jacobian rank 2 from the LMFDB's database of genus 2 curves \cite{bianchi2022rational,lmfdb}.

\subsection{Aim of this paper}
 After Bianchi and Padurariu's computations, the beta version of the LMFDB added a collection of modular curves satisfying the following \cite{lmfdbbeta}:
 
 \begin{LMFDBCondition*}\label{LMFDB Conditions}The LMFDB beta database contains data on modular curves $X_H$ with $\det(H)$ equal to $\hat{\mathbb{Z}}^\times$ such that at least one of the following holds:
    \begin{enumerate}
    \item $H$ has level at most 70.
    \item $H$ has prime-power level at most 335.
    \item $H$ contains $-I$ and has genus at most 24 and level at most 335.
    \item $H$ does not contain $-I$ and has genus at most 8 and level at most 335.
\end{enumerate}
\end{LMFDBCondition*}

The database includes information on rational cusps and rational CM points for every $X_H$. Rational and low degree points data are available for points that correspond to elliptic curves already in the LMFDB, but these points are not yet provably computed for all curves in the database.

In this manuscript, we consider genus 2 bielliptic curves from this database with Jacobian rank 2 as well as curves of Jacobian rank 1 that achieve Jacobian rank 2 over a quadratic imaginary number field with an eye towards understanding the set of \emph{mock rational points}, which we describe below.

 Recall that Balakrishnan and Dogra \cite{Balakrishnan2018QC1} define a finite set $X(\mathbb{Q}_p)_Z\subset\mathbb{Q}_p$, dependent on a choice of correspondence $Z\subset X\times X$, that contains $X(\mathbb{Q}_p)_2$. Our work will use the simplification developed in \cite{bianchi2022rational}: the set of rational points is contained in the finite set 
 \begin{equation}\label{Adef}
     A=\{z\in X(\mathbb{Q}_p):\widetilde{\rho}(z)\in\Omega\}.
 \end{equation}
 
 Here, $\widetilde{\rho}$ is an explicitly computable locally analytic function on $X(\mathbb{Q}_p)$. We will define it precisely in the next section. Additionally, we can explicitly describe and compute $\Omega$ as a finite subset of $\mathbb{Q}_p$ (see Theorem \ref{bianchi2.3}). In this paper, we will refer to the points in $A\setminus X(\mathbb{Q})$ as mock rational points. We aim to study this set with a focus on the question:

\begin{question}
    Are there algebraic points among the mock rational points for a given curve?
\end{question}

The question of algebraic points on modular curves has been studied in recent years by {\"O}zman, Box, and others (see \cite{ozman2019quadratic},\cite{ozman2023quadratic},\cite{box2020quadratic},\cite{box2023cubic}). We have also seen algebraic points appear as mock rational points before. For example, when computing the set of points $X_0(37)(\mathbb{Q}(i))$, Balakrishnan, Dogra, and M\"uller also discovered points on $X_0(37)$ defined over $\mathbb{Q}(\sqrt{-3})$ \cite{Balakrishnan2018QC1}. This example invites several questions. Under what conditions do these points emerge? Over what fields do these points live? What characteristics do these points share?

To explore these questions, we carry out computations on the genus 2 bielliptic modular curves described above. We analyze genus 2 curves with non-simple Jacobian of rank 1 and 2: 209 curves with Jacobian rank 2 and 1237 curves with Jacobian rank 1.

For the genus 2 curves of Jacobian rank 2 over $\mathbb{Q}$, we can use quadratic Chabauty. On the other hand, for curves of Jacobian rank 1 over $\mathbb{Q}$, we consider them over quadratic imaginary fields over which the rank of the Jacobian of the curve increases to 2 and carry out quadratic Chabauty with respect to that field. In particular, since the Jacobian of a bielliptic curve is isogenous to a product of elliptic curves $E_1\times E_2$, the rank is the sum of the ranks of the elliptic curves. We focus on the case where both $E_1$ and $E_2$ have rank 1 over $\mathbb{Q}(\sqrt{-D})$.

In the rank 2 case, we vary the prime used for quadratic Chabauty to observe differences in what points appear. We take $p$ a prime of good, ordinary reduction for $X$ with $3<p<100$. In the rank 1 case, we vary not only the prime but also the quadratic imaginary field used. We determine over which fields $\mathbb{Q}(\sqrt{-D})$ with $1<D<20$ and $D\in\{43,67,163\}$ the rank of the Jacobian of the curve is 2 (while we sought a tractable bound for $D$, adding the three larger values allows us to consider \textit{all} imaginary quadratic fields of class number one at minimal additional cost). For each such field $\mathbb{Q}(\sqrt{-D})$, we restrict to primes $p$ of good, ordinary reduction that split in $\mathbb{Q}(\sqrt{-D})$.

 We provide code and data on GitHub to demonstrate computations in the paper \cite{finnertycode}. The code is modified from \cite{balakrishnancode}, which in turn is modified from \cite{bianchisage} and uses SageMath \cite{sagemath}. We will describe the new modifications in Sections 3 and 4. From these computations, we obtain the following:

\begin{proposition}
    Among the genus 2 Jacobian rank 2 modular curves satisfying \ConditionName{} (\ref{LMFDB Conditions}), Tables \ref{tab:r2exc} and 6\footnote{Tables 6 and 7 are on the Supplementary Tables file on GitHub \cite{finnertycode}.} contain curves that have points over the listed number fields through quadratic Chabauty for some prime $p<100$ of good, ordinary reduction, along with the smallest prime for which they appear, up to precision $O(p^{25})$.
\end{proposition}

\begin{proposition}
    Among the genus 2 Jacobian rank 1 modular curves satisfying \ConditionName{} (\ref{LMFDB Conditions}), Tables \ref{tab:r1exc} and 7 contain curves that have points over the listed number fields through quadratic Chabauty for some prime $p<100$ of good, ordinary reduction, along with the smallest prime for which they appear, up to precision $O(p^{25})$.
\end{proposition}

We note that this list of algebraic points is not provably exhaustive. We can, however, make a conjecture based on the results of the explicit computations. It holds for the genus 2 modular curves satisfying \ConditionName{} (Section \ref{LMFDB Conditions}).

\begin{conjecture}\label{conj1}
    Let $X$ be a bielliptic modular curve of genus 2 and level $N$. If $Jac(X)$ is of rank 2 over $\mathbb{Q}$ and $X$ has a mock rational point over $\mathbb{Q}(\sqrt{D})$ then $D$ must be a divisor of $N$.
\end{conjecture}

We might hope for an analogous conjecture in the rank one case. Suppose $X$ is a bielliptic modular curve of genus 2 and level $N$ with $Jac(X)$ of rank 1. If $X$ has a mock rational point over $\mathbb{Q}(\sqrt{D})$ arise through quadratic Chabauty when computing over $\mathbb{Q}(\sqrt{D'})$, we frequently observe one of the following:
    \begin{enumerate}
        \item $D$ divides $N$.
        \item $D$ divides $D'$.
    \end{enumerate}
    Note that in order to use $\mathbb{Q}(\sqrt{D'})$ for quadratic Chabauty contributions for the curve $X$, we need the rank of the Jacobian of $X$ to be 2 over this field.

However, we have observed examples for which this is not the case. The curve \curvepage{24.36.2.bo.1} has $\mathbb{Q}(\sqrt{-3})$-points appear when using $p=13$ and working over $\mathbb{Q}(\sqrt{-10})$. Also, the curve $X_0(37)$ has the aforementioned points over $\mathbb{Q}(\sqrt{-3})$ appear when using $p=13$ and working over $\mathbb{Q}(i)$. We suspect there is a systematic explanation for these not captured in the above divisibility conditions on $D$.

\subsection{Acknowledgements}
Many thanks are due to Francesca Bianchi and Oana Padurariu for their virtual correspondence related to their previous work. We also thank Jennifer Balakrishnan for proposing the project and for helpful and insightful conversations throughout. The third modular curves workshop for the Simons Collaboration in Arithmetic Geometry, Number Theory, and Computation in March 2024, with support from the Simons Foundation, was also instrumental. Lastly, we are grateful to the reviewers of this article for their insight and attention to detail, both of which greatly improved the paper. The author was partially supported by NSF Grant DMS-1945452 and US-Israel BSF Grant 2022393.

\section{Quadratic Chabauty for Bielliptic Curves}
Our main source for this expository section is \cite{bianchi2022rational}. Let $X$ be a non-singular genus 2 bielliptic curve 
    \begin{equation}\label{biellform}
        X:y^2=a_6x^6+a_4x^4+a_2x^2+a_0,a_i\in\mathbb{Z},
    \end{equation}
and consider the two elliptic curves
   \begin{equation}\label{E1form}
        E_1:y^2=x^3+a_4x^2+a_2a_6x+a_0a_6^2, 
    \end{equation}
    \begin{equation}\label{E2form}
        E_2:y^2=x^3+a_2x^2+a_4a_0x+a_6a_0^2.
    \end{equation}
We have degree 2 maps $\phi_i:X\rightarrow E_i$ with
\begin{equation}
    \phi_1(x,y)=(a_6x^2,a_6y), \hspace{0.5cm}\phi_2(x,y)=(a_0x^{-2},a_0yx^{-3}).
\end{equation}

\begin{remark}
    The genus 2 curves of our dataset do not initially come in this form. However, any genus 2 curve with a degree 2 map to an elliptic curve has a model of this form \cite{wetherell97thesis}.
\end{remark}

Recall our goal is to define a function $\widetilde{\rho}$ so that we can compute the set $A$ in Equation \ref{Adef}. Let $\log$ denote the $p$-adic logarithm extended to $\mathbb{Q}_p^\times$ by setting $\log(p)=0$. While some summands of $\widetilde{\rho}$ will depend on the choice of branch of $\log$, the final function is independent of the value of $\log(p)$. On $E_i$ (although we could replace this discussion with any $E/\mathbb{Q}$) with good reduction at $p$, for $P\in E_i(\mathbb{Q}_p)$, we let $Log(P)=\int_\infty^P\omega_i$, where $\omega_i=\frac{dx}{2y}$ and the integral is Coleman integration \cite[Theorem 2.8]{coleman1985torsion}.

If $P_1$ and $P_2$ are in the same residue disk, then $Log(P_1)-Log(P_2)=\int_{P_1}^{P_2}\omega_i$ and the integral can be calculated via formal antidifferentiation. Since $\omega_i$ is holomorphic, the function $Log$ is locally analytic. Furthermore, $Log$ vanishes at a point $P$ if and only if that point is torsion \cite[Proposition 3.1]{coleman1985torsion}.

We consider the decomposition of the global $p$-adic height into local heights $\lambda_q$. Given a prime $q$, we let $\text{ord}_q$ denote the $q$-adic valuation on $\mathbb{Q}_q$, normalized to be surjective onto the integers and let $|\cdot|_q$ denote the standard absolute value on $\mathbb{Q}_q$. With this, the following results give us the desired objects:

\begin{proposition}[\hspace{-1pt}{\cite[Proposition 2.2]{bianchi2022rational}}]
    Let $q\neq p$. If $P\in E(\mathbb{Q}_q)$ reduces to a nonsingular point modulo $q$, then $\lambda_q(P)=\log(max\{1,|x(P)|_q\})$. 

    Let $W_q^E$ be the set of values attained by $\lambda_q$ on points in $E(\mathbb{Q}_q)$ of the form $(x,y)$ with $x,y\in\mathbb{Z}_q$. Then $W_q^E$ is finite, explicitly computable, and equal to $\{0\}$ for all but finitely many $q$. 

    In particular, $W_q^E\subseteq\{0\}$ at all primes of good reduction for the given model of the elliptic curve $E$.
\end{proposition}

Note that in the above proposition, we require the model for $E$ to be minimal, see \cite{bianchi2020bielliptic}. Next, we use the following additional notation in the following theorem and remark. For a prime $q$, we define
\begin{equation*}
    Z_q = X(\mathbb{Q}_q)\setminus\{P:x(P)\in\{0,\infty\}\}.
\end{equation*}

\begin{theorem}[\hspace{-1pt}{\cite[Theorem 2.3]{bianchi2022rational}}\label{bianchi2.3}]
    Suppose that each of $E_1$ and $E_2$ have rank 1 over $\mathbb{Q}$, and let $p$ be a prime of good reduction for the equation of $X$. For each $i\in\{1,2\}$, fix a choice of subspace of $H^1_{dR}(E_i/\mathbb{Q}_p)$ complementary to the space of holomorphic forms, and consider the corresponding global height $h_p$ and local N\'eron functions $\lambda_q$ at every $q$. Let $P_i\in E_i(\mathbb{Q})$ be a point of infinite order. Define $\alpha_i=\frac{h_p(P_i)}{Log^2(P_i)}$. Then:
    \begin{enumerate}
        \item The constant $\alpha_i$ is independent of choice of $P_i$.
        \item The function $\rho:Z_p\rightarrow\mathbb{Q}_p$ given by
        \begin{equation*}
            \rho(z)=\lambda_p(\phi_1(z))-\lambda_p(\phi_2(z))-2\log(x(z))-\alpha_1Log^2(\phi_1(z))+\alpha_2Log^2(\phi_2(z))
        \end{equation*}
        can be continued to a locally analytic function $\widetilde{\rho}:X(\mathbb{Q}_p)\rightarrow\mathbb{Q}_p$.
        \item For a prime $q\neq p$, let
        \begin{equation*}
        \begin{aligned}
            \Omega_q=(-W_q^{E_1}+W_q^{E_2}+\{-n\log q:-\text{ord}_q(a_6)\leq n\leq \text{ord}_q(a_0), n\equiv0\bmod2\}&) \\ \cup (\log|a_0|_q-W_q^{E_1})\cup(-\log|a_6|_q+W_q^{E_2}&) 
        \end{aligned} 
        \end{equation*}
        and set $\Omega=\{\sum_{q\text{ bad}}w_q:w_q\in\Omega_q\}$. Then $\Omega$ is finite and contains the image $\widetilde{\rho}(X(\mathbb{Q}))$.
    \end{enumerate}
\end{theorem}

\begin{remark}
    In the development of the theory and the first examples of quadratic Chabauty for bielliptic curves, Balakrishnan and Dogra also define a set $\Omega$ \cite[Theorem 1.2]{Balakrishnan2018QC1}. In this paper, the set of points whose image under a certain quadratic Chabauty function lands in $\Omega$ is a superset rather than strictly equal to $X(\mathbb{Q}_p)_2$ (note that Balakrishnan and Dogra worked over general number fields). In turn, the set of Theorem \ref{bianchi2.3} is a superset of the set of Balakrishnan and Dogra, so the set $A$ of Equation \ref{Adef} is a superset of $X(\mathbb{Q}_p)_2$. This method of giving a superset that is explicitly computable is standard: in addition to \cite{Balakrishnan2018QC1}, see \cite{balakrishnan2019cursed}, \cite{balakrishnan2020QCII}. The authors of \cite{balakrishnan2019cursed} remark that this can be thought of as an analogue to working with a suitable quotient of the Jacobian as in \cite{mazur1977} and \cite{merel1996}.
\end{remark} 

\section{Computations over the Rationals}

We begin with 209 curves of Jacobian rank 2. Many have the same model, and we use the coarse curves only to obtain a unique list of models. This strategy is permissible because any fine modular curve in our list shares a model with a coarse modular curve. This fact can be checked by hand: for example, the fine curve \curvepage{14.96.2-14.g.1.1} has Weierstrass model 
\begin{equation*}
    y^2 = x^6 - 3x^5 - x^4 + 7x^3 - x^2 - 3x + 1,
\end{equation*} which is the same as that of the coarse curve \curvepage{14.48.2.g.1}. By excluding fine modular curves, we obtain a list of 76 equations defining bielliptic curves. We also exclude curves with no known rational points. We then note that many modular curves have exceptional isomorphisms. We consider only one representative from each isomorphism class, giving us 30 equations $y^2=f(x)$ to consider \footnote{While we checked for these isomorphisms by noting when Weierstrass models on the LMFDB were identical or off by a sign difference, checks for exceptional isomorphisms were also performed by David Roe in October 2024 \cite{roe2024isomorphisms}, and our results agree for the curves under consideration here.}. We convert each to a bielliptic model by using Magma to obtain Richelot isogenous surfaces that are of the form $E_1\times E_2$. Performing changes of variables on these elliptic curves allows us to find $a_i$ that determine a model for $X$ in the form of Equation \ref{biellform}. We run the analysis for each prime $p$ less than or equal to 97 of good, ordinary reduction.

\subsection{On Choice of Bielliptic Model}
A given curve does not have a unique representation of the form in Equation \ref{biellform}. Most immediately, a change of variables $[x:y:z]\mapsto[-z:y:x]$ gives an isomorphism between
    $y^2=a_6x^6+a_4x^4+a_2x^2+a_0$
and
    $y^2=a_0x^6+a_2x^4+a_4x^2+a_6$.
However, there are often other possibilities. For example, the following are isomorphic:
\begin{align*}
    y^2&=-x^6 + 83x^4 - 19x^2 + 1, \\
    y^2&=64x^6 - 304x^4 + 332x^2 - 1.
\end{align*}
A model with a larger discriminant will produce a larger set of $p$-adic points when completing quadratic Chabauty. We are therefore interested in finding a model of minimal discriminant to streamline the analysis.

The following result describes how discriminants of isomorphic curves are related to each other. For more, one can refer to \cite{liu1994}.
\begin{Lemma}[\hspace{-1pt}{\cite[Section 2.1]{Booker_Sijsling_Sutherland_Voight_Yasaki_2016}}]
    Let $F$ be a perfect field and $X$ a curve of genus 2 over $F$. Suppose $X$ has a Weierstrass equation $y^2+h(x)y=f(x)$. If $(y')^2+h'(x')y'=f'(x')$ is another Weierstrass model for the curve $X$, then it is related to the original via a change of variables of the form
    \begin{equation*}
        x'=\frac{ax+b}{cx+d},\hspace{1cm}y'=\frac{ey+j(x)}{(cx+d)^3},\hspace{1cm}\text{with }ad-bc,e\in F^\times,j\in F[x],\deg j\leq3.
    \end{equation*}
    Furthermore, we have the discriminants $\Delta(f',h')$ and $\Delta(f,h)$ are related as follows:
    \begin{equation*}
        \Delta(f',h')=e^{20}(ad-bc)^{-30}\Delta(f,h).
    \end{equation*}
\end{Lemma}

The last statement of the lemma immediately gives us the following:

\begin{corollary}[\hspace{-1pt}{\cite[Section 2.1]{Booker_Sijsling_Sutherland_Voight_Yasaki_2016}}]\label{discinv}
    The discriminant is invariant as an element of $F^\times/F^{\times10}$.
\end{corollary}

We now describe how to attain a model of minimal discriminant in our case. We can represent transformations between models with matrices of the form
\begin{equation*}
    M(a,b,c,d,e):=\begin{pmatrix}
    a & 0 & b\\0&e&0\\c&0&d\end{pmatrix}.
\end{equation*}

By the lemma, this transformation changes the discriminant by a factor of $e^{20}(ad-bc)^{-30}$. The following is a result of direct computation with the previous lemma and the definition of $M(a,b,c,d,e)$.

\begin{proposition}\label{transformationcondition}
    Let $X:y^2=a_6x^6+a_4x^4+a_2x^2+a_0$ with $a_i$ all integers. For $M(a,b,c,d,e)$ to yield an isomorphic curve of the form in Equation \ref{biellform} with integer coefficients, it is necessary and sufficient to have the following conditions:
    \begin{itemize}
        \item The following are all equal to zero:
        \tiny{\begin{itemize}
       \item $a_0(3ac^5) + a_2(2ac^3d^2 + bc^4d) + a_4(acd^4 + 2bc^2d) + a_6(3bd^5)$.
        \item $a_0(5a^3c^3) + a_2(a^3cd + 3a^2bc^2d + ab^2c^3) + a_4(a^2bd^3 +  3ab^2cd^2 + b^3c^2) + a_6(5b^3d^3)$.
        \item $a_0(3a^5c) + a_2(a^4bd + 2a^3b^2c) + a_4(2a^2b^3d + ab^4c) + a_6(3b^5d)$.
          
        \end{itemize}}
        \normalsize{\item The following are all congruent to zero modulo $(ad-bc)^6$:}
       \tiny{ \begin{itemize}
            \item $e^2(c^6a_0+c^4d^2a_2 + c^2d^4a_4 + d^6a_6)$
            \item $ e^2(15a^2c^4a_0 + 6a^2c^2d^2a_2 + a^2d^4a_4 + 8abc^3da_2 + 8abcd^3a_4 + b^2c^4a_2 + 6b^2c^2d^2a_4 + 15b^2d^4a_6)$
            \item $e^2(15a^4c^2a_0 + a^4d^2a_2 + 8a^3bcda_2 + 6a^2b^2c^2a_2 + 6a^2b^2d^2a4 + 8ab^3cda_4 + b^4c^2a_4 + 15b^4d^2a_6)$
            \item $e^2(a^6a_0 + a^4b^2a_2 + a^2b^4a_4 + b^6a_6)$
        \end{itemize}}
    \end{itemize}
\end{proposition}

We can think of the first set as vanishing conditions: these ensure that the odd-power terms have zero coefficient. The latter are integrality conditions that ensure the even-power terms are elements of $\mathbb{Z}$.

Let $\prod_ip_i^{e_i}$ be the prime factorization of the discriminant of a given model for the curve $X$. We compute the potential discriminant-lowering effects $p^k<1$ and, for each such $p^k$, the set of $M(a,b,c,d,e)$ such that $e^{20}(ad-bc)^{-30}=p^k$ subject to the condition that the absolute values of $a,b,c,d,e$ are less than or equal to the maximum $p_i$. In addition to considering integral $a,b,c,d,e$, we also consider numbers of the form $1/p$ for $p$ prime to ensure that we can obtain transformations with effects that reduce the discriminant by a prime power that is a negative even multiple of 10. 

Note that our bounds on the matrix coefficients are sufficient because we can obtain a reduction in the discriminant by a factor of $p^{10}$ for any $p$ dividing the discriminant, which is sufficient by Corollary \ref{discinv}. While matrices that reduce the discriminant by a given factor are not unique, we offer the following as an example of a transformation to reduce the discriminant by $p^{-10}$:
\begin{equation*}
    M(p,0,0,1,p)=\begin{pmatrix}
        p & 0 & 0 \\ 0 & p & 0 \\0 & 0 & 1
\end{pmatrix}.
\end{equation*}

For each potential factor by which we might lower the discriminant, we apply Proposition \ref{transformationcondition} to obtain a matrix $M(a,b,c,d,e)$ that will yield an isomorphic model of the desired form with integer coefficients, if such a matrix exists. If it exists, we apply it and obtain a new model. If not, we continue to the next factor. We repeat until there is no more allowable transformation matrix that would reduce the discriminant while keeping it in the desired form. In this way, we obtain optimal models to use in the analysis.

For example, initial transformation of the curve \curvepage{30.60.2.i.1} into a bielliptic model yielded
\begin{equation*}
    y^2=-2700x^6 + 8000,
\end{equation*}
which has discriminant $2^{54}3^{21}5^{25}$. Applying the algorithm produces the model
\begin{equation*}
    y^2=-108x^6+5,
\end{equation*}
of discriminant $2^{24}3^{21}5^5$, removing a factor of $2^{30}5^{20}$.

\subsection{Execution of the Analysis}

A description of the algorithm is as follows. Our input is: $f$ a polynomial over $\mathbb{Q}$ of the form in Equation \ref{biellform} such that the two elliptic curves of Equations \ref{E1form} and \ref{E2form} both have rank 1 over $\mathbb{Q}$, $p$ an odd prime such that $E_1$ and $E_2$ have good reduction and a fixed $p$-adic precision $n$. In our computations, we set $n=25$, sufficient for all curves under the precision analysis in \cite{bianchi2022rational}. Our output is a list of rational points on $H: y^2=f(x)$ and a list of $p$-adic points on $H$ that are not recognized as rational points. We compute this by the following:
\begin{enumerate}
    \item For each elliptic curve $E_i$, compute the following: $E_{i,m}$ (the minimal models), $\Delta_i$ (the discriminant), $h_i$ (the canonical $p$-adic height), $P_i$ (a generator of the free part of $E_i(\mathbb{Q})$), $E_{i,K}$ (the elliptic curve over $K=\mathbb{Q}_p$ up to precision $n$), $E_{i,p}$ (the elliptic curve over $\mathbb{F}_p$). 
    \item Set $H_p$ as $H:y^2=f(x)$ over $\mathbb{F}_p$. Let $D$ be the set of points of $H_p$ up to automorphism.
    \item Compute the set $\Omega$ per Theorem \ref{bianchi2.3}. 
    \item For every point in $D$, compute $\widetilde{\rho}$ per Theorem \ref{bianchi2.3}. Then, for each element $\omega\in\Omega$, compute the zeroes of $\widetilde{\rho}-\omega$. Note that these will be $p$-adic numbers. Find the points on $H$ over $\mathbb{Q}_p$ with these zeroes as $x$-coordinates and save them in a list $L$. 
    \item For each element of $L$, attempt to lift to a rational point. Return a list of rational points and a list of candidate algebraic points for further analysis.
\end{enumerate}

We focus on $p$-adic points that do not obviously lift to rational points. We use SageMath's algdep functionality to see if either of the points' coordinates satisfy  relations of degree up to four. We distinguish between Weierstrass points, which are unsurprising to encounter, and other non-Weierstrass points.

\begin{remark}
    Since we work with precision $n=25$, one might wonder whether the points truly appear via this analysis or if they are only transcendental points that are $p$-adically very close. In all cases, we can immediately confirm that the algebraic point satisfies the defining equation for the curve. Then, to be more robust, we can apply the Weierstrass preparation theorem. We consider the power series given by $\widetilde\rho-\omega$ and factor out the polynomial part that has a root corresponding to the point in question. If we can evenly factor out that polynomial part with coefficients in $\mathbb{Q}$ up to our working precision and the root is simple, then we have the desired result.

    For example, consider the curve \curvepage{16.48.2.a.1}, denoted $X$, with model given by $y^2=x^6 - 5x^4 - 5x^2 + 1$ of Jacobian rank 2. Working with both $p=11$ and $p=17$, we find $p$-adic points that appear to correspond to $(\pm1,\pm2\sqrt{-2})$. Without loss of generality, we let $Q=(1,2\sqrt{-2})$, since the other three points are redundant up to automorphism. One can immediately see that $Q$ is in $X(\mathbb{Q}(\sqrt{-2}))$ (as are its images under automorphism).

    Let $t_0$ represent the root of $\widetilde\rho-\omega$ that was recognized as $Q$ for $p=17$. It turns out that $(t-t_0)$ evenly divides $\widetilde\rho-\omega$. This is a simple root with a linear coefficient of just 1, so we can conclude that the point is algebraic. Although the terms in the case of $p=11$ are 11-adic instead of 17-adic, the result of being able to evenly factor out a simple linear term is the same.
\end{remark}

Another change to the code relates to the computation of the set $\Omega$. When computing the local height at a given bad prime $q$, one can check if the curve has potential good reduction at $q$. If so, then the contribution at that prime can be understood entirely by information about the local height at a single $\mathbb{Q}_q$-point on the curve \cite{bianchi2022rational}. The original implementation assumes that such a point will exist. We remove this assumption by first checking in MAGMA \cite{MR1484478} if the curve is locally soluble at $q$ . In the event that the code fails to obtain a $\mathbb{Q}_q$-point, we proceed to the general case to compute the height contribution at $q$. 

The last change is concerning the point at infinity. At each residue disc, we need to compute the function $\widetilde{\rho}$. For the given point $P$, which is a $\mathbb{F}_p$-point, we let $Q$ denote a lift to the curve over $\mathbb{Q}_p$. We let $f_1$ and $f_2$ denote the images of $Q$ to $E_1$ and $E_2$ over $\mathbb{Q}_p$. If, for one $f_i$, we have $f_i$ multiplied by the order of $E_i$ over $\mathbb{F}_p$ is the point at infinity, then we must compute the Coleman integral from the point at infinity to $f_i$. This step can be prohibitively computationally expensive for high precision. Therefore, we take advantage of the fact that $\int_{\infty}^{P_1}\omega_i=\frac{1}{2}\int_{-P_1}^{P_1}\omega_i$ to circumvent difficulties with the point at infinity. This is a consequence of extending the following lemma to tangential basepoints, as was done on page 289 of \cite{balakrishnan2011appendix}:
\begin{Lemma}[\hspace{-1pt}{\cite[Lemma 16]{balakrishnan2010coleman}}]
  Let $\omega_i$ be an odd, everywhere meromorphic differential on $X$. Choose $P,Q$ as points in $ X(\mathbb{C}_p)$ which are not poles of $\omega_i$, with $P$ Weierstrass. Then for $\iota$ the hyperelliptic involution, we have $\int_P^Q\omega_i=\frac{1}{2}\int_{\iota(Q)}^Q\omega_i$.  
\end{Lemma}

After the analysis described in this section, we were also able to replicate the strategy of Bianchi and Padurariu to use the Mordell-Weil sieve to confirm the complete sets of rational points for these curves. For information about the methodology, one can refer to Sections 4.3 and 4.4 of \cite{bianchi2022rational}. This allowed us to compute the rational points on more curves, but since the technology behind the computation is not new, it is not a focus of this paper. It does, however, give the following:

\begin{proposition}
    The analysis of \cite{finnertycode} produces the set of rational points for genus 2 bielliptic modular curves of Jacobian rank 1 or 2 that satisfy \ConditionName{} (Section \ref{LMFDB Conditions}).
\end{proposition}

\subsection{Algebraic Weierstrass Mock Rational Points}
Weierstrass points form a convenient set of examples when considering in which cases an algebraic point of a curve will appear in the analysis. We obtain precise criteria by comparing the splitting behavior of the defining polynomials with the results in Table 6, which can be found in the Supplementary Tables file of the GitHub \cite{finnertycode}.

\begin{proposition}\label{r2Weierstrass}
    Let $P=(x_0,0)$ be a Weierstrass point of a genus 2 bielliptic modular curve $X$ with rank 2 Jacobian satisfying \ConditionName{} (\ref{LMFDB Conditions}) with $x_0\notin\mathbb{Q}$. Then for all $p < 100$, all nonrational Weierstrass points $P$ in $A\setminus X(\mathbb{Q})$ (see Equation \ref{Adef}), up to precision $O(p^{25})$, satisfy the following properties:
    \begin{enumerate}
        \item The prime $p$ is of good and ordinary reduction for both elliptic curves $E_i$, where $Jac(X)\sim E_1\times E_2$, so it satisfies the conditions for quadratic Chabauty for bielliptic curves.
        \item The defining polynomial $f(x)$ for $X$ factors as $(ax^2+bx+c)g(x)$, and $x_0$ is a root of $ax^2+bx+c$.
        \item The point $P$ lives in $\mathbb{Q}(\sqrt{D})$ and $D$ is a quadratic residue modulo $p$. 
    \end{enumerate}
\end{proposition}

Consider the curve $X$ with label \curvepage{24.72.2.iy.1}. It has a model given by $y^2=6x^6 + 12x^4 + 12x^2 + 6$. The defining polynomial for $X$ factors over the rationals as $6 (x^2 + 1) (x^2 - x + 1) (x^2 + x + 1)$, so we have three eligible terms in the factorization, and we expect to observe all 6 Weierstrass points. Note $E$ does not have ordinary reduction at 7 or 47, so we cannot use these primes for quadratic Chabauty. All other primes $p$ with $3<p<100$ are of good, ordinary reduction. Since the Weierstrass points live over $\mathbb{Q}(i)$, and $\mathbb{Q}(\sqrt{-3})$, we can check the Legendre symbol of these primes with respect to $-1$ and $-3$. The behavior we see is as expected: the points $(\pm i,0)$ appear exactly when $p\equiv1\bmod4$ and the points $(\pm\frac{1}{2}\pm\frac{\sqrt{-3}}{2},0)$ appear exactly when $p\equiv1\bmod 3$.

Next, consider the modular curve \curvepage{16.48.2.a.1}. We use $y^2=x^6-5x^5-5x^2+1$, which factors over the rationals as $(x^2+1)(x^2-2x-1)(x^2+2x-1)$. The first term has roots $\pm i$, and the latter terms have roots living over $\mathbb{Q}(\sqrt{2})$. The primes $p$ that are 1 modulo 4 are those such that $-1$ is a quadratic residue modulo $p$, so we consider the list $\{5, 13, 17, 29, 37, 41, 53, 61, 73, 89, 97\}$. We can check which primes are of good, ordinary reduction for the corresponding elliptic curves and obtain the list $\{17, 41, 73, 89, 97\}$. These are exactly the primes for which we observe the points $(\pm i,0)$. However, we do not observe the points over $\mathbb{Q}(\sqrt{2})$ for any prime $<100$. The conditions for a Weierstrass point to appear are necessary but not sufficient.

\section{Computations over Quadratic Number Fields}
We perform a similar process as the start of the analysis of curves of Jacobian rank 2. Although we start with 1237 curves with Jacobian rank 1, we obtain a list of 116 unique defining equations of curves with a known rational point. We once again convert them to the desired bielliptic form $y^2=a_6x^6+a_4x^4+a_2x^2+a_0$. In this section, we use code from Bianchi and Padurariu \cite{padurariumagma} to assist in some of the transformations. 

As stated in the introduction, an eligible choice of number field $K=\mathbb{Q}(\sqrt{-D})$ and prime $p$ is one where the Jacobian of the curve is rank 2 over $K$ and $p$ is a prime of good, ordinary reduction for the curve that splits in $K$. For each $p$ between 3 and 100 of good, ordinary reduction, we check potential values of $D$ until we have obtained a list of at least five eligible pairs or run out of primes. We obtain a list of eligible pairs $(D,p)$ for each curve and run quadratic Chabauty on each.

Code modifications from the previous section were also used for this analysis. In addition, we need to consider the height contributions of primes that divide the discriminant of the number field $K$. We incorporate this into the computation of $\Omega$. More precisely, we can modify Theorem \ref{bianchi2.3} as follows:

\begin{theorem}\label{rank1rhoomega}
        Let $X$ be a bielliptic genus 2 curve of the form in Equation \ref{biellform} whose Jacobian splits into a product of two elliptic curves as $J\sim E_1\times E_2$ such that $rk(E_1(\mathbb{Q}))+rk(E_2(\mathbb{Q}))=1$. 
        
        Suppose that each of $E_1$ and $E_2$ has rank 1 over an imaginary quadratic number field $K/\mathbb{Q}$, and let $p$ be a prime of good reduction for the equation of $X$. For each $i\in\{1,2\}$, fix a choice of subspace of $H^1_{dR}(E_i/\mathbb{Q}_p)$ complementary to the space of holomorphic forms, and consider the corresponding global height $h_p$ and local N\'eron functions $\lambda_q$ at every $q$. 
        
        Let $P_i\in E_i(K)$ be a point of infinite order. Define $\alpha_i=\frac{h_p(P_i)}{[K:\mathbb{Q}]Log^2(P_i)}$. Then:
    \begin{enumerate}
        \item The constant $\alpha_i$ is independent of choice of $P_i$.
        \item The function $\rho:Z_p\rightarrow\mathbb{Q}_p$ given by
        \begin{equation*}
            \rho(z)=\lambda_p(\phi_1(z))-\lambda_p(\phi_2(z))-2\log(x(z))-\alpha_1Log^2(\phi_1(z))+\alpha_2Log^2(\phi_2(z))
        \end{equation*}
        can be continued to a locally analytic function $\widetilde{\rho}:X(\mathbb{Q}_p)\rightarrow\mathbb{Q}_p$.
        \item For a prime $q\neq p$, let
        \begin{equation*}
        \begin{aligned}
            \Omega_q=(-W_q^{E_1}+W_q^{E_2}+\{-n\log q:-\text{ord}_q(a_6)\leq n\leq \text{ord}_q(a_0), n\equiv0\bmod2\}&) \\ \cup (\log|a_0|_q-W_q^{E_1})\cup(-\log|a_6|_q+W_q^{E_2}&)
        \end{aligned} 
        \end{equation*}
        and set $\Omega=\{\sum_{q\in S}w_q:w_q\in\Omega_q\}$, where $S$ consists of the bad primes for each $E_i$, the prime factors of $a_0$ and $a_6$ for $X$, \textbf{and} the prime divisors of the discriminant of $K$. Then $\Omega$ is finite and contains $\widetilde{\rho}(X(\mathbb{Q}))$.
    \end{enumerate}
\end{theorem}

Many of the arguments of Bianchi and Padurariu for Theorem \ref{bianchi2.3} in the case of $K=\mathbb{Q}$ apply, see \cite{bianchi2022rational}. We now discuss the necessary steps to make this extended result hold. The addition of the factor $\frac{1}{[K:\mathbb{Q}]}$ in the definition of $\alpha_i$ is a necessary normalization factor for $p$-adic heights over a quadratic field, see Corollary 8.1 ii) of \cite{Balakrishnan2018QC1}. With this normalization complete, the only difference with the previous Theorem is the elements of the set $S$ that contribute to the definition of $\Omega$. In particular, the set $S$ encapsulates the possible nontrivial local height contributions. When working over a number field $K$, it is possible that the ramification of primes can lead to nontrivial height contributions, meaning we need to consider them. However, since there are only finitely many such primes, $\Omega$ still remains finite. 

\subsection{Appearance of Weierstrass Points in the Rank 1 Case}
    We checked the rank 1 results to compare the appearance of Weierstrass points to that in the rank 2 case, keeping in mind the results described in Proposition \ref{r2Weierstrass}. Since we did not run the analysis on Jacobian rank 1 curves for all primes between 3 and 100, we cannot directly claim an analogue to Proposition \ref{r2Weierstrass}. However, for primes that were tested, Weierstrass points that appeared satisfied the conditions, regardless of the number field used.

    Consider the curve \curvepage{60.48.2.c.1} with model $y^2=x^6 + 20x^4 + 25x^2 - 750$ as an example. This factors as $-2 (x^2 - 5) (x^2 + 10) (x^2 + 15)$. The Weierstrass points live over $\mathbb{Q}(\sqrt{5}),\mathbb{Q}(\sqrt{-10})$, and $\mathbb{Q}\sqrt{-15})$. The primes used were 11, 13, 19, 23, and 37. As expected, we saw the points over $\mathbb{Q}(\sqrt{5})$ when using $p=11$ and 19, the points over $\mathbb{Q}(\sqrt{-10})$ using all five primes, and the points over $\mathbb{Q}(\sqrt{-15})$ using $p=19$ and 23.

\section{Results}

Table \ref{tab:r2exc} lists those curves of Jacobian rank 2 over $\mathbb{Q}$ for which the analysis yielded non-Weierstrass algebraic \textbf{non-rational} points. For each curve, the table states:
\begin{enumerate}
    \item The identifying label of the curve, linked to its page on the Beta LMFDB.
    \item The model used in the analysis.
    \item The smallest prime that realized the points.
    \item The number field over which the points live.
    \item The set of discovered algebraic points with respect to the model used.
\end{enumerate}Table 6 enumerates the curves of Jacobian rank 2 with detected Weierstrass points. Table \ref{tab:r1exc} and Table 7 do the same for the Jacobian rank 1 curves. However, instead of simply stating the minimal prime, we also state the value $D$ for which the computation was performed over $\mathbb{Q}(\sqrt{-D})$. Table \ref{Summary} summarizes the results. Note that the sum of each row is greater than the number of unique models because some curves had points in multiple categories. Tables \ref{tab:r2exc} and \ref{tab:r1exc} are in the appendix and Tables 6 and 7 are on the Supplementary Tables file on GitHub \cite{finnertycode}.

As discussed in the computations sections, some modular curves with different labels share a model. We remark on labeling conventions here \cite{lmfdb}. For coarse modular curves $X_H$, labels have the form N.i.g.c.n, which denote respectively the level, index, genus, Gassmann class identifier, and an integer to distinguish nonconjugate subgroups where the first four identifiers are the same. Table 8 in the Supplementary Tables on GitHub lists each label whose model realized quadratic points with any other modular curves in the database that share a model \cite{finnertycode}. This includes both exceptional isomorphisms and fine modular curves. 

Any such curve has level divisible by that of the curve used in the main results. For example, our first curve is \curvepage{14.48.2.g.1}. In our analysis, we used the model $y^2=-x^6+83x^4-19x^2+1$. This model can also represent several other curves, including \curvepage{28.96.2-14.g.1.1}, \curvepage{42.96.2-14.g.1.1}, \curvepage{56.96.2-14.g.1.1}, and \curvepage{70.96.2-14.g.1.1}. We can add these as evidence for Conjecture \ref{conj1}.

\begin{table}[H]\caption{Count of genus 2 curves for which a point with the given characteristic was found, by rank 
 of the Jacobian}\label{Summary}
\begin{tabular}{l|l|l|l|l}
Rank &
  \begin{tabular}[c]{@{}l@{}}$x$-coordinate\\ irrational, $y\neq0$\end{tabular} &
  \begin{tabular}[c]{@{}l@{}}$x$-coordinate rational, \\ $y$-coordinate irrational\end{tabular} &
  Weierstrass &
  None \\ \hline
2 &
  6 &
  3 &
  14 &
  10 \\
1 &
  42 &
  8 &
  39 &
  65
\end{tabular}
\end{table}

\section{Examples}
\subsection{$X_{ns}^+(15)$}
We highlight points found on the modular curve $X_{ns}^+(15)$, which has LMFDB label \curvepage{15.60.2.d.1}. The model we use is $y^2=-45x^6 + 75x^4 - 15x^2 + 1$. The primes less than 100 of good, ordinary reduction for the curve are 7, 13, 19, 31, 37, 43, 61, 67, 73, and 97. For most of these primes, we recover no notable mock rational points. However, for $p\in\{19,31,61\}$ we observe four points: $\{(\pm\frac{\sqrt{5}}{5},\pm\frac{4}{5})\}$. The sizes of the relevant outputs are in Table \ref{ns15example}.

\begin{table}[b]\caption{Analysis Data for $X_{ns}^+(15)$}\label{ns15example}
\begin{tabular}{l|l|l}
prime & \begin{tabular}[c]{@{}l@{}}\# recovered \\rational points\end{tabular}& \begin{tabular}[c]{@{}l@{}}\# mock rational \\$p$-adic points\end{tabular} \\ \hline
19    & 14                & 284                    \\
31    & 14                & 404                    \\
61    & 14                & 780                   
\end{tabular}
\end{table}

\subsection{A Rank 1 Example}
We highlight points over different number fields found on the modular curve \curvepage{16.24.2.f.1}, using model $y^2 = 4x^6 + 6x^4 + 4x^2 + 1$. The Jacobian of this curve has rank 1. The results for two iterations of the analysis are in Table \ref{16example}. Both cases recovered the two known rational points. Note that the points $(\pm i,\pm i)$ appeared both times, but the points over $\mathbb{Q}(\sqrt{-7})$ only appeared in the latter. A given viable field and prime will not necessarily capture all algebraic points.

    \begin{table}[]\caption{Analysis Data for Curve \curvepage{16.24.2.f.1}}\label{16example}
\begin{tabular}{l|l|l|l}
\multicolumn{1}{c|}{prime} &
  \multicolumn{1}{c|}{\begin{tabular}[c]{@{}c@{}}number \\ field\end{tabular}}  &
  \multicolumn{1}{c|}{\begin{tabular}[c]{@{}c@{}}\# mock rational\\ $p$-adic points\end{tabular}} &
  \multicolumn{1}{c}{Algebraic points} \\ \hline
5 &
  $\mathbb{Q}(\sqrt{-6})$ &
  66 &
  $(\pm i,\pm i)$ \\
37 &
  $\mathbb{Q}(\sqrt{-7})$ &
  566 &
  $(\pm i,\pm i), (\pm\frac{2\sqrt{-7}}{7},\pm\frac{5\sqrt{-7}}{49})$
\end{tabular}
\end{table}

\subsection{A Return to $X_0(37)$}
We conclude this section by revisiting an example from the beginning of this paper. A motivating example for performing the analysis of curves with Jacobian rank 1 was $X_0(37)$, which has label \curvepage{37.38.2.a.1}. Recall it was found to have points over both $\mathbb{Q}(i)$ and $\mathbb{Q}(\sqrt{-3})$. In this work, we found further points over $\mathbb{Q}(\sqrt{-7})$.  

Using the model $y^2 = -x^6-9x^4-11x^2+37$, we can run the analysis with $p=11$ over $\mathbb{Q}(\sqrt{-7})$. In this case, we observe the following set of points:
\begin{equation*}
    \{(\pm\sqrt{-7},\pm4),(\pm\frac{\sqrt{-7}}{3},\pm\frac{172}{27})\}.
\end{equation*}
The interested reader can find more about the quadratic points on $X_0(37)$ in Balakrishnan and Mazur's 2024 paper regarding Ogg's Torsion Conjecture \cite{balakrishnan2024oggs} and an earlier paper of Box \cite{box2020quadratic}.

\section{Future Work}
We have gathered extensive data regarding quadratic points that appear during quadratic Chabauty computations for modular curves. We conclude with some remarks regarding additional analyses and questions to explore. These questions emerged through conversations with Kiran Kedlaya and the anonymous referees. 

First, a potential counterexample to the conjecture could be found in the following setup. Given a rational Heegner point on one of the elliptic curves $E_1$ or $E_2$ in the decomposition of the Jacobian of the genus 2 curve, this could pull back to a quadratic point on $X$. If this point satisfies the local height conditions required to be a mock rational point, it could be possible to appear in the quadratic Chabauty analysis while not having the properties of the points in the conjecture. These examples, if they exist, would be rare. It is worth exploring if they are possible, and to find them if so.

We also have two questions immediately from looking at the tables:

\begin{question}
   All points found were over quadratic or biquadratic fields. Does this generalize to modular curves of different genus or ranks of the Jacobian, or to other hyperelliptic curves in general?
\end{question}
\begin{question}
    It is also the case that upon observing one algebraic point, we always find the other points in its Galois orbit. Once again, how does this generalize?
\end{question}  
Another interesting set of questions question pertains to the Mordell-Weil group of the Jacobian.
\begin{question}
    Could some of the algebraic points found be linearly equivalent to those in the Mordell-Weil group of the Jacobian over the base field? If they have the same heights away from $p$, could this explain their appearance in the analysis?
\end{question}
A related but not entirely equivalent question is the following:
\begin{question}
    If we find algebraic points over number fields for a given curve, how does the Mordell-Weil group change when computed over these fields or their compositum versus the rationals? 
\end{question}

\section{Appendix: Tables}
\small{
\begin{longtable}{|l|l|l|l|l|}
\caption{Rank 2 Non-Weierstrass Points} \label{tab:r2exc} \\
\hline \multicolumn{1}{|c|}{\textbf{Label}}  & \multicolumn{1}{c|}{\textbf{Model $y^2=f(x)$}}& \multicolumn{1}{c|}{\textbf{$p$}} & \multicolumn{1}{c|}{\textbf{Number Field}}& \multicolumn{1}{c|}{\textbf{Points}}\\ \hline 
\endfirsthead

\hline \multicolumn{5}{|r|}{{Continued on next page}} \\ \hline
\endfoot

\hline \hline
\endlastfoot
\curvepage{14.48.2.g.1}   & $-x^6 + 83x^4 - 19x^2 + 1$    & 37    &$\mathbb{Q}(\sqrt{-3},\sqrt{21})$ & $(\pm\frac{\sqrt{-3}}{3},\pm\frac{8\sqrt{21}}{9})$ \\
\curvepage{15.60.2.d.1}   & $-45x^6 + 75x^4 - 15x^2 + 1$ & 19    &$\mathbb{Q}(\sqrt{5})$& $(\frac{\sqrt{5}}{5},\pm\frac{4}{5})$               \\
\curvepage{16.48.2.a.1}  & $x^6 - 5x^4 - 5x^2 + 1$ &11    &$\mathbb{Q}(\sqrt{-2})$& $(\pm1,\pm2\sqrt{2})$       \\
\curvepage{18.72.2.f.1} & $9x^6 - 99x^4 + 27x^2 - 1$   & 7     &$\mathbb{Q}(\sqrt{-3}  )$& $(\pm\frac{\sqrt{-3}}{3},\pm\frac{8\sqrt{-3}}{3})$          \\
\curvepage{28.48.2.j.1} &$x^6 - 83x^4 + 19x^2 - 1$  & 19    &$\mathbb{Q}(\sqrt{-3},\sqrt{7})$ & $(\pm\frac{\sqrt{-3}}{3},\pm\frac{8\sqrt{-21}}{9})$       \\
\curvepage{36.72.2.d.1}&$-9x^6 + 99x^4 - 27x^2 + 1$   & 13    &$\mathbb{Q}(\sqrt{-3},\sqrt{3}  )$  & $(\pm\frac{\sqrt{-3}}{3},\pm\frac{8\sqrt{3}}{3})$\\
\curvepage{48.48.2.dk.1} &$12x^6 + 54x^4 + 48x^2 + 12$& 73    &$\mathbb{Q}(\sqrt{3} )$  & $(0,\pm2\sqrt{3})$  \\
\curvepage{48.48.2.dm.1}&$24x^6 + 54x^4 - 24x^2 + 3$ & 73    &$\mathbb{Q}(\sqrt{3})$ &$(0,\pm\sqrt{3})$       \\
\curvepage{48.48.2.dw.1} &$27x^6 - 72x^4 + 54x^2 - 8$& 73    &$\mathbb{Q}(\sqrt{6}  )$&$(\pm\frac{\sqrt{6}}{9},\pm2)$
\end{longtable}
}

\scriptsize{
\begin{center}
\begin{longtable}{|l|l|l|l|l|}
\caption{Rank 1 Non-Weierstrass Points} \label{tab:r1exc} \\

\hline \multicolumn{1}{|c|}{\textbf{Label}}&\multicolumn{1}{c|}{\textbf{Model $y^2=f(x)$}}  & \multicolumn{1}{c|}{\textbf{$(D,p)$}} & \multicolumn{1}{c|}{\textbf{Number Field}}&\multicolumn{1}{c|}{\textbf{Points}}\\ \hline 
\endfirsthead

\multicolumn{5}{c}
{{\bfseries \tablename\ \thetable{} -- continued from previous page}} \\
\hline \multicolumn{1}{|c|}{\textbf{Label}}&\multicolumn{1}{c|}{\textbf{Model}}  & \multicolumn{1}{c|}{\textbf{$(D,p)$}} & \multicolumn{1}{c|}{\textbf{Number Field}}&\multicolumn{1}{c|}{\textbf{Points}}\\ \hline 
\endhead

\hline \multicolumn{5 }{|r|}{{Continued on next page}} \\ \hline
\endfoot

\hline \hline
\endlastfoot
\curvepage{11.66.2.a.1}  & $11x^6 + 11x^4 - 7x^2 + 1$ & [7,23]  &$\mathbb{Q}(\sqrt{-7})$  & $(\pm\frac{\sqrt{-7}}{3},\frac{76}{27})$      \\
\curvepage{16.24.2.f.1}& $4x^6 + 6x^4 + 4x^2 + 1$  & [6,5]   &$\mathbb{Q}(i)$& $(\pm i,\pm i)$                 \\
\curvepage{16.24.2.f.1}  & $4x^6 + 6x^4 + 4x^2 + 1$ & [7,37]  &$\mathbb{Q}(\sqrt{-7})$  & $(\pm\frac{2\sqrt{-7}}{7},\pm\frac{5\sqrt{-7}}{49})$      \\
\curvepage{16.48.2.bm.1} &$8x^6 + 16x^4 + 9x^2 + 1$ & [5,41]  &$\mathbb{Q}(\sqrt{-5})$          &$(\pm\frac{3\sqrt{-5}}{7},\pm\frac{62}{343})$\\
\curvepage{16.48.2.bx.1} &$-x^6 + 7x^4 - 7x^2 + 1$& [1,5]   &$\mathbb{Q}(i)$ &$(\pm i,\pm4)$                 \\
\curvepage{16.48.2.bx.1} &$-x^6 + 7x^4 - 7x^2 + 1$& [17,7] &$\mathbb{Q}(\sqrt{-17})$&$(\pm\sqrt{-17},\pm84)$ \\
\curvepage{16.48.2.c.1} & $-x^6 - 5x^4 + 5x^2 + 1$ & [7,11]  &$\mathbb{Q}(\sqrt{-7})$ &$(\pm\sqrt{-7},\pm8)$          \\
\curvepage{16.48.2.c.1}  & $-x^6 - 5x^4 + 5x^2 + 1$ & [1,17]  &$\mathbb{Q}(i,\sqrt{2})$  &$(\pm i,\pm2\sqrt{-2})$      \\
\curvepage{16.48.2.c.2} &$-2x^6 - 10x^4 + 10x^2 + 2$ & [19,11] &$\mathbb{Q}(\sqrt{-19})$   &$(\pm\frac{9\sqrt{-19}}{19},\pm\frac{680\sqrt{-19}}{361})$    \\
\curvepage{16.48.2.c.2}  & $-2x^6 - 10x^4 + 10x^2 + 2$&[1,17]  &$\mathbb{Q}(i)$   &$(\pm i,\pm 4i)$               \\
\curvepage{16.48.2.y.1} &$-8x^6 + 18x^4 - 8x^2 + 1$  & [15,17] &$\mathbb{Q}(\sqrt{2})$ &$(\pm\frac{\sqrt{2}}{2},\pm\frac{\sqrt{2}}{2})$         \\
\curvepage{18.54.2.e.1} &$-3x^6 + 18x^4 - 3x^2 + 4$ & [3,7]   &$\mathbb{Q}(\sqrt{-3})$ &\makecell[l]{$(\pm\frac{\sqrt{-3}}{3},\pm\frac{8}{3}),$\\$(\pm\sqrt{-3},\pm16)$ }   \\
\curvepage{20.30.2.l.1} &$256x^6 + 113x^4 + 18x^2 + 1$  & [7,11]  &$\mathbb{Q}(\sqrt{-7})$    &$(\pm\frac{\sqrt{-7}}{7},\pm\frac{2\sqrt{-7}}{49})$       \\
\curvepage{24.36.2.a.1}& $8x^6 + 1$  & [7,37]  &$\mathbb{Q}(\sqrt{-7})$  &$(\pm\frac{2\sqrt{-7}}{7},\pm\frac{13\sqrt{-7}}{49})$    \\
\curvepage{24.36.2.bo.1}&$72x^6 + 48x^4 + 12x^2 + 1$ & [10,13] &$\mathbb{Q}(\sqrt{-3})$ &$(\pm\frac{\sqrt{-3}}{3},\pm\frac{\sqrt{-3}}{3})$          \\
\curvepage{24.36.2.bo.1} &$72x^6 + 48x^4 + 12x^2 + 1$ & [2,19]  &$\mathbb{Q}(\sqrt{-2})$ &$(\pm\frac{\sqrt{-2}}{2},\pm\sqrt{-2})$         \\
\curvepage{24.36.2.bv.1} &$-2x^6 + 2$& [1,13]  &$\mathbb{Q}(i)$   &$(\pm i,\pm2)$               \\
\curvepage{24.36.2.ev.1} &$x^6 - 13x^4 + 39x^2 - 27$& [3,13]  &$\mathbb{Q}(\sqrt{-3})$   &$(0,\pm3\sqrt{-3})$    \\
\curvepage{24.36.2.ev.1} &$x^6 - 13x^4 + 39x^2 - 27$& [3,19]  &$\mathbb{Q}(\sqrt{-3},\sqrt{-2})$&$(\pm\sqrt{-3},\pm12\sqrt{-2})$ \\
\curvepage{24.36.2.gc.1}&$-6x^6 + 6$ & [1,13]  &$\mathbb{Q}(\sqrt{3},i)$&$(\pm i,\pm2\sqrt{3})$\\
\curvepage{24.36.2.gh.1} &$-6x^6 - 78x^4 + 78x^2 + 6$& [1,5]   &$\mathbb{Q}(i)$  & $(\pm i,\pm12i)$             \\
\curvepage{24.36.2.gj.1} &$-9x^6 - 21x^4 + 21x^2 + 9$& [1,5]   &$\mathbb{Q}(i,\sqrt{-6})$ &$(\pm i,\pm2\sqrt{-6})$        \\
\curvepage{24.36.2.gk.1} &$x^6 + 13x^4 - 13x^2 - 1$ & [1,5]   &$\mathbb{Q}(i,\sqrt{6})$&$(\pm i,\pm2\sqrt{6})$        \\
\curvepage{24.36.2.gk.1} &$x^6 + 13x^4 - 13x^2 - 1$& [1,13]  &$\mathbb{Q}(i)$   & $(0,\pm i)$              \\
\curvepage{24.36.2.gm.1} &$6x^6 + 14x^4 - 14x^2 - 6$& [1,5]   &$\mathbb{Q}(i)$  &$(\pm i,\pm 4)$           \\
\curvepage{24.36.2.gm.1}  &$6x^6 + 14x^4 - 14x^2 - 6$ &[11,5]   &$\mathbb{Q}(\sqrt{-11})$  &$(\pm\frac{\sqrt{-11}}{5},\pm\frac{192}{125})$              \\
\curvepage{24.36.2.x.1} &$2x^6 - 2$& [1,13]   &$\mathbb{Q}(i)$  &$(\pm i,\pm2i)$                 \\
\curvepage{24.36.2.x.1} &$2x^6 - 2$& [17,7]   &$\mathbb{Q}(\sqrt{-17})$  &$(\pm\frac{15\sqrt{-17}}{17},\pm\frac{4774\sqrt{-17}}{289})$               \\
\curvepage{24.72.2.ba.1} &$24x^6 + 12x^4 + 18x^2 + 9$& [8,43]  &$\mathbb{Q}(\sqrt{-2})$ &$(\pm\frac{\sqrt{-2}}{4},\pm\frac{21}{8})$          \\
\curvepage{24.72.2.cm.1} &$-9x^6 + 3x^4 - 3x^2 + 1$& [1,13]  &$\mathbb{Q}(i)$     &$(\pm i,\pm4)$     \\
\curvepage{24.72.2.da.1} &$-72x^6 + 12x^4 - 6x^2 + 1$& [8,43]  &$\mathbb{Q}(\sqrt{-2})$ &$(\pm\frac{\sqrt{-2}}{2},\pm4)$          \\
\curvepage{24.72.2.db.1} &$9x^6 - 6x^4 + 12x^2 - 8$& [10,13] &$\mathbb{Q}(\sqrt{-10})$&$(\pm\frac{\sqrt{-10}}{15},\pm\frac{208\sqrt{-10}}{225})$ \\
\curvepage{24.72.2.db.1} &$9x^6 - 6x^4 + 12x^2 - 8$& [2,19]  &$\mathbb{Q}(\sqrt{-2})$&\makecell[l]{$(0,\pm2\sqrt{-2}),$\\$(\pm\sqrt{-2},\pm8\sqrt{-2})$}    \\
\curvepage{24.72.2.fy.1} &$216x^6 + 72x^4 + 12x^2 + 1$& [17,13] &$\mathbb{Q}(\sqrt{-3})$  &$(\pm\frac{\sqrt{-3}}{3},\pm\sqrt{-3})$         \\
\curvepage{24.72.2.hy.1} &$-6x^6 + 12x^4 - 12x^2 + 6$& [1,5]   &$\mathbb{Q}(i)$   &$(\pm i,\pm6)$                \\
\curvepage{24.72.2.ja.1} &$24x^6 - 24x^4 + 12x^2 - 3$& [2,11]  &$\mathbb{Q}(\sqrt{-2})$        &\makecell[l]{$(\pm\frac{\sqrt{-2}}{2},\pm3\sqrt{-2})$\\$(\pm\frac{1}{2},\pm\frac{3\sqrt{-2}}{4})$}  \\
\curvepage{24.72.2.t.1}  & $-24x^6 + 12x^4 - 18x^2 + 9$&[6,31]  &$\mathbb{Q}(\sqrt{-6})$ &$(\pm\frac{\sqrt{-6}}{2},\pm12)$\\
\curvepage{24.72.2.u.1}  &$9x^6 + 9x^4 + 3x^2 + 3$& [1,13]  &$\mathbb{Q}(i)$  &$(\pm2i,\pm21i)$                 \\
\curvepage{28.32.2.a.1} &$2x^6 + 46x^4 + 30x^2 - 14$& [19,5]  &$\mathbb{Q}(\sqrt{-19})$  &$(\pm\sqrt{-19},\pm48)$        \\
\curvepage{28.32.2.a.1}   &$2x^6 + 46x^4 + 30x^2 - 14$& [2,17]  &$\mathbb{Q}(\sqrt{-2})$          &$(\pm\frac{3\sqrt{-2}}{2},\pm\frac{49}{2})$ \\
\curvepage{28.48.2.g.1} &$49x^6 + 70x^4 + 25x^2 + 4$ & [17,13] &$\mathbb{Q}(\sqrt{-17})$  &$(\pm\frac{\sqrt{-17}}{7},\pm\frac{64}{49})$  \\
\curvepage{28.48.2.g.1}  &$49x^6 + 70x^4 + 25x^2 + 4$ & [67,19] &$\mathbb{Q}(\sqrt{7})$        &$(\pm\frac{\sqrt{7}}{7},\pm\frac{8\sqrt{7}}{7})$   \\
\curvepage{30.30.2.b.1} &$150x^6 - 690x^4 + 546x^2 - 6$ & [19,17] &$\mathbb{Q}(\sqrt{-15})$  & $(\pm\frac{\sqrt{-15}}{5},\pm\frac{32\sqrt{-15}}{5})$\\
\curvepage{36.24.2.a.1} &$x^6 + 108$ & [15,19] &$\mathbb{Q}(\sqrt{-3})$ &$(\pm\sqrt{-3},\pm9)$      \\
\curvepage{36.54.2.i.1} &$x^6 + 6x^4 + 9x^2 + 16$ & [3,13]  &$\mathbb{Q}(\sqrt{-3})$       &$(\pm\sqrt{-3},\pm4)$  \\
\curvepage{37.38.2.a.1}&$-x^6 - 9x^4 - 11x^2 + 37$ & [3,7]   &$\mathbb{Q}(\sqrt{-3})$       &\makecell[l]{$(\pm\sqrt{-3},\pm4)$,\\$(\pm3\sqrt{-3},\pm116)$}  \\
\curvepage{37.38.2.a.1}  &$-x^6 - 9x^4 - 11x^2 + 37$ & [7,11]  &$\mathbb{Q}(\sqrt{-7})$       &\makecell[l]{$(\pm\sqrt{-7},\pm4$)\\$(\pm\frac{\sqrt{-7}}{3},\pm\frac{172}{27})$}\\
\curvepage{38.40.2.b.1} &$4x^6 + 64x^4 + 304x^2 + 361$ & [3,13]  &$\mathbb{Q}(\sqrt{-3})$  &$(\pm\frac{19\sqrt{-3}}{12},\pm\frac{437\sqrt{-3}}{288})$         \\
\curvepage{40.30.2.e.1} &$25x^6 - 115x^4 + 91x^2 - 1$ & [1,13]  &$\mathbb{Q}(i)$   &$(0,\pm i)$              \\
\curvepage{40.30.2.f.1} & $x^6 - 91x^4 + 115x^2 - 25$ & [1,13]  &$\mathbb{Q}(i)$          &$(0,\pm5i)$ \\
\curvepage{45.54.2.c.1} & $5x^6 + 21x^4 - 9x^2 - 1$ & [1,13]  &$\mathbb{Q}(i)$  & $(0,\pm i)$   \\
\curvepage{48.48.2.ez.1} & $3x^6 - 21x^4 + 21x^2 - 3$ & [1,13] & $\mathbb{Q}(\sqrt{-3},i)$ &$(\pm i,\pm4\sqrt{-3})$ \\
\curvepage{56.48.2.l.1}&$2x^6 + 25x^4 + 140x^2 + 196$ & [17,13] &$\mathbb{Q}(\sqrt{14})$&$(\pm\sqrt{14},\pm112)$\\
\curvepage{60.30.2.p.1} &$-3x^6 - 5x^4 - 225x^2 - 3375$ & [17,7] &$\mathbb{Q}(\sqrt{-5})$&$(\pm\sqrt{-5},\pm20\sqrt{-5})$ \\
\curvepage{60.30.2.p.1} &$-3x^6 - 5x^4 - 225x^2 - 3375$ & [15,19] &$\mathbb{Q}(\sqrt{-15})$&$(\pm\frac{\sqrt{-15}}{6},\pm\frac{355\sqrt{-15}}{24})$ \\
\curvepage{60.48.2.f.1} &$9x^6 + 60x^4 + 25x^2 - 250$ & [13,17] &$\mathbb{Q}(\sqrt{-15})$&$(\pm\frac{\sqrt{-15}}{3},\pm\frac{10\sqrt{-15}}{3})$ \\
\curvepage{70.30.2.b.1} &$x^6 - 98x^4 + 833x^2 - 10976$ & [7,11]  &$\mathbb{Q}(\sqrt{-7})$  &$(\pm\sqrt{-7},\pm56\sqrt{-7})$         
\end{longtable}

\end{center}}
\normalsize{
\bibliographystyle{alphaurl}
\bibliography{sources}{}

\newcommand{\etalchar}[1]{$^{#1}$}
\begin{thebibliography}{AKMJ{\etalchar{+}}24}

\bibitem[AKMJ{\etalchar{+}}24]{ozman2023quadratic}
Nikola Ad\v{z}aga, Timo Keller, Philippe Michaud-Jacobs, Filip Najman, Ekin Ozman, and Borna Vukorepa.
\newblock Computing quadratic points on modular curves {$X_0(N)$}.
\newblock {\em Math. Comp.}, 93(347):1371--1397, 2024.
\newblock \href {https://doi.org/10.1090/mcom/3902} {\path{doi:10.1090/mcom/3902}}.

\bibitem[Bal]{balakrishnancode}
Jennifer Balakrishnan.
\newblock Sagemath code.
\newblock URL: \url{https://github.com/jbalakrishnan/QC_bielliptic}.

\bibitem[BBBM21]{balakrishnan2021numberfields}
Jennifer~S. Balakrishnan, Amnon Besser, Francesca Bianchi, and J.~Steffen M\"uller.
\newblock Explicit quadratic {C}habauty over number fields.
\newblock {\em Israel J. Math.}, 243(1):185--232, 2021.
\newblock \href {https://doi.org/10.1007/s11856-021-2158-5} {\path{doi:10.1007/s11856-021-2158-5}}.

\bibitem[BBK10]{balakrishnan2010coleman}
Jennifer~S. Balakrishnan, Robert~W. Bradshaw, and Kiran~S. Kedlaya.
\newblock Explicit {C}oleman integration for hyperelliptic curves.
\newblock In {\em Algorithmic number theory}, volume 6197 of {\em Lecture Notes in Comput. Sci.}, pages 16--31. Springer, Berlin, 2010.
\newblock URL: \url{https://doi.org/10.1007/978-3-642-14518-6_6}, \href {https://doi.org/10.1007/978-3-642-14518-6\_6} {\path{doi:10.1007/978-3-642-14518-6\_6}}.

\bibitem[BCP97]{MR1484478}
Wieb Bosma, John Cannon, and Catherine Playoust.
\newblock The {M}agma algebra system. {I}. {T}he user language.
\newblock {\em J. Symbolic Comput.}, 24(3-4):235--265, 1997.
\newblock Computational algebra and number theory (London, 1993).
\newblock URL: \url{http://dx.doi.org/10.1006/jsco.1996.0125}, \href {https://doi.org/10.1006/jsco.1996.0125} {\path{doi:10.1006/jsco.1996.0125}}.

\bibitem[BD18]{Balakrishnan2018QC1}
Jennifer~S. Balakrishnan and Netan Dogra.
\newblock Quadratic {C}habauty and rational points, {I}: {$p$}-adic heights.
\newblock {\em Duke Math. J.}, 167(11):1981--2038, 2018.
\newblock With an appendix by J. Steffen M\"uller.
\newblock \href {https://doi.org/10.1215/00127094-2018-0013} {\path{doi:10.1215/00127094-2018-0013}}.

\bibitem[BD21]{balakrishnan2020QCII}
Jennifer~S. Balakrishnan and Netan Dogra.
\newblock Quadratic {C}habauty and rational points {II}: {G}eneralised height functions on {S}elmer varieties.
\newblock {\em Int. Math. Res. Not. IMRN}, (15):11923--12008, 2021.
\newblock \href {https://doi.org/10.1093/imrn/rnz362} {\path{doi:10.1093/imrn/rnz362}}.

\bibitem[BDM{\etalchar{+}}19]{balakrishnan2019cursed}
Jennifer Balakrishnan, Netan Dogra, J.~Steffen M\"uller, Jan Tuitman, and Jan Vonk.
\newblock Explicit {C}habauty-{K}im for the split {C}artan modular curve of level 13.
\newblock {\em Ann. of Math. (2)}, 189(3):885--944, 2019.
\newblock \href {https://doi.org/10.4007/annals.2019.189.3.6} {\path{doi:10.4007/annals.2019.189.3.6}}.

\bibitem[BGG23]{box2023cubic}
Josha Box, Stevan Gajovi\'c, and Pip Goodman.
\newblock Cubic and quartic points on modular curves using generalised symmetric {C}habauty.
\newblock {\em Int. Math. Res. Not. IMRN}, (7):5604--5659, 2023.
\newblock \href {https://doi.org/10.1093/imrn/rnab358} {\path{doi:10.1093/imrn/rnab358}}.

\bibitem[Bia20]{bianchi2020bielliptic}
Francesca Bianchi.
\newblock Quadratic {C}habauty for (bi)elliptic curves and {K}im's conjecture.
\newblock {\em Algebra Number Theory}, 14(9):2369--2416, 2020.
\newblock \href {https://doi.org/10.2140/ant.2020.14.2369} {\path{doi:10.2140/ant.2020.14.2369}}.

\bibitem[BKK11]{balakrishnan2011appendix}
Jennifer~S. Balakrishnan, Kiran~S. Kedlaya, and Minhyong Kim.
\newblock Appendix and erratum to ``{M}assey products for elliptic curves of rank 1'' [mr2629986].
\newblock {\em J. Amer. Math. Soc.}, 24(1):281--291, 2011.
\newblock \href {https://doi.org/10.1090/S0894-0347-2010-00675-3} {\path{doi:10.1090/S0894-0347-2010-00675-3}}.

\bibitem[BM25]{balakrishnan2024oggs}
Jennifer~S. Balakrishnan and Barry Mazur.
\newblock Ogg's torsion conjecture: fifty years later.
\newblock {\em Bull. Amer. Math. Soc. (N.S.)}, 62(2):235--268, 2025.
\newblock With an appendix by Netan Dogra.
\newblock \href {https://doi.org/10.1090/bull/1851} {\path{doi:10.1090/bull/1851}}.

\bibitem[Box21]{box2020quadratic}
Josha Box.
\newblock Quadratic points on modular curves with infinite {M}ordell-{W}eil group.
\newblock {\em Math. Comp.}, 90(327):321--343, 2021.
\newblock \href {https://doi.org/10.1090/mcom/3547} {\path{doi:10.1090/mcom/3547}}.

\bibitem[BPa]{padurariumagma}
Francesca Bianchi and Oana Padurariu.
\newblock {MAGMA} code.
\newblock URL: \url{https://github.com/oana-adascalitei/MWSieveForDatabase}.

\bibitem[BPb]{bianchisage}
Francesca Bianchi and Oana Padurariu.
\newblock Sage{M}ath code.
\newblock URL: \url{https://github.com/bianchifrancesca/QC_bielliptic}.

\bibitem[BP24]{bianchi2022rational}
Francesca Bianchi and Oana Padurariu.
\newblock Rational points on rank 2 genus 2 bielliptic curves in the {LMFDB}.
\newblock In {\em Lu{C}a{NT}: {LMFDB}, computation, and number theory}, volume 796 of {\em Contemp. Math.}, pages 215--242. Amer. Math. Soc., Providence, RI, 2024.
\newblock \href {https://doi.org/10.1090/conm/796/16003} {\path{doi:10.1090/conm/796/16003}}.

\bibitem[BSS{\etalchar{+}}16]{Booker_Sijsling_Sutherland_Voight_Yasaki_2016}
Andrew~R. Booker, Jeroen Sijsling, Andrew~V. Sutherland, John Voight, and Dan Yasaki.
\newblock A database of genus-2 curves over the rational numbers.
\newblock {\em LMS J. Comput. Math.}, 19:235--254, 2016.
\newblock \href {https://doi.org/10.1112/S146115701600019X} {\path{doi:10.1112/S146115701600019X}}.

\bibitem[Cha41]{chabauty1941points}
Claude Chabauty.
\newblock Sur les points rationnels des courbes alg\'ebriques de genre sup\'erieur \`a{} l'unit\'e.
\newblock {\em C. R. Acad. Sci. Paris}, 212:882--885, 1941.

\bibitem[Cola]{lmfdb}
The~LMFDB Collaboration.
\newblock The {L}-functions and modular forms database.
\newblock URL: \url{http://www.lmfdb.org}.

\bibitem[Colb]{lmfdbbeta}
The~LMFDB Collaboration.
\newblock The {L}-functions and modular forms database, {B}eta version.
\newblock URL: \url{http://beta.lmfdb.org}.

\bibitem[Col85a]{coleman1985effective}
Robert~F. Coleman.
\newblock Effective {C}habauty.
\newblock {\em Duke Math. J.}, 52(3):765--770, 1985.
\newblock \href {https://doi.org/10.1215/S0012-7094-85-05240-8} {\path{doi:10.1215/S0012-7094-85-05240-8}}.

\bibitem[Col85b]{coleman1985torsion}
Robert~F. Coleman.
\newblock Torsion points on curves and {$p$}-adic abelian integrals.
\newblock {\em Ann. of Math. (2)}, 121(1):111--168, 1985.
\newblock \href {https://doi.org/10.2307/1971194} {\path{doi:10.2307/1971194}}.

\bibitem[Fal83]{faltings1983proof}
G.~Faltings.
\newblock Endlichkeitss\"atze f\"ur abelsche {V}ariet\"aten \"uber {Z}ahlk\"orpern.
\newblock {\em Invent. Math.}, 73(3):349--366, 1983.
\newblock \href {https://doi.org/10.1007/BF01388432} {\path{doi:10.1007/BF01388432}}.

\bibitem[Fin25]{finnertycode}
Kate Finnerty.
\newblock Sage{M}ath code, 2025.
\newblock URL: \url{https://github.com/kfinnerty19/QC_bielliptic}.

\bibitem[Kim09]{kim2009Selmer}
Minhyong Kim.
\newblock The unipotent {A}lbanese map and {S}elmer varieties for curves.
\newblock {\em Publ. Res. Inst. Math. Sci.}, 45(1):89--133, 2009.
\newblock \href {https://doi.org/10.2977/prims/1234361156} {\path{doi:10.2977/prims/1234361156}}.

\bibitem[Liu94]{liu1994}
Qing Liu.
\newblock Conducteur et discriminant minimal de courbes de genre {$2$}.
\newblock {\em Compositio Math.}, 94(1):51--79, 1994.
\newblock URL: \url{http://www.numdam.org/item?id=CM_1994__94_1_51_0}.

\bibitem[Maz77]{mazur1977}
B.~Mazur.
\newblock Modular curves and the {E}isenstein ideal.
\newblock {\em Inst. Hautes \'Etudes Sci. Publ. Math.}, (47):33--186, 1977.
\newblock With an appendix by Mazur and M. Rapoport.
\newblock URL: \url{http://www.numdam.org/item?id=PMIHES_1977__47__33_0}.

\bibitem[Mer96]{merel1996}
Lo\"ic Merel.
\newblock Bornes pour la torsion des courbes elliptiques sur les corps de nombres.
\newblock {\em Invent. Math.}, 124(1-3):437--449, 1996.
\newblock \href {https://doi.org/10.1007/s002220050059} {\path{doi:10.1007/s002220050059}}.

\bibitem[Mor23]{mordell1922rational}
L.~J. Mordell.
\newblock On the rational solutions of the indeterminate equations of the third and fourth degrees.
\newblock {\em Proc. Cambridge Philos. Soc.}, 21:179--192, 1922/23.

\bibitem[OS19]{ozman2019quadratic}
Ekin Ozman and Samir Siksek.
\newblock Quadratic points on modular curves.
\newblock {\em Math. Comp.}, 88(319):2461--2484, 2019.
\newblock \href {https://doi.org/10.1090/mcom/3407} {\path{doi:10.1090/mcom/3407}}.

\bibitem[Roe24]{roe2024isomorphisms}
David Roe.
\newblock Private communication.
\newblock Harvard University CMSA Math and Machine Learning Program Zulip Server, 2024.

\bibitem[{The}24]{sagemath}
{The Sage Developers}.
\newblock {\em {S}ageMath, the {S}age {M}athematics {S}oftware {S}ystem ({V}ersion 9.3)}, 2024.
\newblock {\tt https://www.sagemath.org}.

\bibitem[Wet97]{wetherell97thesis}
Joseph~Loebach Wetherell.
\newblock {\em Bounding the number of rational points on certain curves of high rank}.
\newblock ProQuest LLC, Ann Arbor, MI, 1997.
\newblock Thesis (Ph.D.)--University of California, Berkeley.
\newblock URL: \url{http://gateway.proquest.com/openurl?url_ver=Z39.88-2004&rft_val_fmt=info:ofi/fmt:kev:mtx:dissertation&res_dat=xri:pqdiss&rft_dat=xri:pqdiss:9803394}.

\end{thebibliography}
}
\end{document}